
\baselineskip=14pt
\parskip=10pt

\font\eightrm=cmr8 

\magnification=\magstephalf

\def\N{{\cal N}}

\def\1{{\overline{1}}}
\def\2{{\overline{2}}}
\parindent=0pt
\overfullrule=0in

\def\frac#1#2{{#1 \over #2}}
\centerline
{\bf 
Experimenting with Standard Young Tableaux
}
\bigskip
\centerline
{\it Shalosh B. EKHAD and Doron ZEILBERGER}
\centerline
\qquad \qquad \qquad 
{\it In fond memory of Herbert Wilf and Albert Nijenhuis, and in honor of Curtis Greene}
\bigskip

{\bf Abstract}:  
Using Symbolic Computation with Maple, we can discover lots of (rigorously-proved!) facts about Standard Young Tableaux, in particular the distribution of the entries in any specific cell, and the sorting probabilities.

{\bf Maple package}

This article is accompanied by a Maple package, {\tt SYT.txt}, available from:

{\tt https://sites.math.rutgers.edu/\~{}zeilberg/tokhniot/SYT.txt} \quad .

The web-page of this article 

{\tt https://sites.math.rutgers.edu/\~{}zeilberg/mamarim/mamarimhtml/syt.html} \quad,

contains many input and output files, some of which will be referred to later.

{\bf Preface} 

One of the most {\it iconic} objects in mathematics, both {\it concrete} [K], and {\it abstract} [F], are Standard Young Tableaux [Wi].
Recall that an {\it integer partition}, or {\it partition}, for short, aka {\it shape}, of a non-negative integer $n$, with $k$ {\bf parts}, is a non-increasing sequence of {\bf positive} integers

$$
\lambda=(\lambda_1 ,\lambda_2,\dots, \lambda_k) \quad,
$$

$$
\lambda_1 \geq \lambda_2 \geq \dots \geq \lambda_k >0 \quad,
$$

such that
$$
\lambda_1 + \dots + \lambda_k \, = \, n \quad.
$$

The {\bf Ferrers diagram} (or {\bf Young diagram})  of a partition $\lambda$ is a left-justified array of dots (or empty boxes) where the top row has $\lambda_1$ dots (boxes),
the second row has $\lambda_2$ boxes, $\dots$, the $k$-th row has $\lambda_k$ dots. For example, the Ferrers diagram of $(4,4,3,1)$ is
$$
\matrix{
*&*&*&*\cr
*&*&*&*\cr
*&*&*\cr
*\cr
} \quad .
$$

Given a shape $\lambda=(\lambda_1, \dots, \lambda_k)$ with $n$ boxes, a {\bf standard Young  tableau} is a way of filling the boxes with the integers $\{1, \dots, n\}$, such that each of
them shows up (necessarily  once) and both rows and columns are increasing. More formally, it is an array
$$
T_{i,j} \quad, \quad 1 \leq i \leq k \quad, \quad 1 \leq j \leq \lambda_i \quad ,
$$
such that $T_{i,j}<T_{i,j+1}$ and $T_{i,j}<T_{i+1,j}$ whenever they exit.

Here are the  five standard Young tableaux of shape $(2,2,1)$:
$$
\matrix{
1&2\cr
3&4\cr
5} 
\quad, \quad
\matrix{
1&2\cr
3&5\cr
4} 
\quad, \quad
\matrix{
1&3\cr
2&4\cr
5} 
\quad, \quad
\matrix{
1&3\cr
2&5\cr
4} 
\quad, \quad
\matrix{
1&4\cr
2&5\cr
3} 
\quad .
$$

{\eightrm To see the set of standard Young tableaux of shape {\tt L} in our Maple package, type {\tt SYT(L);} . For example to see the above five tableaux type {\tt SYT([2,2,1]);} .}

The total number of standard Young tableaux of shape $\lambda$, denoted by $f^{\lambda}$,  is famously given by the {\it hook length formula}, or equivalently (and more convenient for us) by the
{\bf Young-Frobenius formula} (see [K]).

$$
f^{\lambda} \, = \,
\frac{(\lambda_1+ \lambda_2+ \dots + \lambda_k)!}{(\lambda_1+k-1)!(\lambda_2+k-2)! \cdots \lambda_k!} \cdot
\prod_{i=1}^k \prod_{j=i+1}^k (\lambda_i-\lambda_j+j-i) \quad.
$$

Fix a shape $\lambda$ and fix a cell $[i,j]$, $1 \leq i \leq k, 1 \leq j \leq \lambda_i$.

Who can be the occupant of that cell? 

Calling that occupant $r$, we have:

$$
ij \, \leq \, r \, \leq \,
\lambda_1+\lambda_2+ \dots + \lambda_{i-1}
+\lambda'_1+\lambda'_2+ \dots + \lambda'_{j-1} -(i-1)(j-1)+1
 \quad ,
$$
where $\lambda'=(\lambda'_1, \dots , \lambda'_{\lambda_1})$ is the {\bf conjugate partition}.

This gives a certain {\bf probability distribution}. What is it?

For example, with the shape $\lambda=(2,2,1)$ again, and the cell $[2,1]$, the set of possible occupants is $\{2,3\}$, and the probability of it being $2$ is $\frac{3}{5}$ and of it being $3$ is $\frac{2}{5}$.

Later on we will be interested not in {\it specific} shapes, but in general (mostly rectangular) shapes, with a fixed number of rows, but arbitrary (i.e. {\it symbolic}) shape.
Fixing the number of rows to be $k$ (where $k$ is {\it numeric}), and regarding the shape $(n,n,\dots, n)$, where $n$ is repeated $k$ times,
we would be interested in deriving closed-form expressions (as rational functions of $n$), for the probability distribution of the possible occupants of a given
first row cell $[1,j]$, for any given {\bf numeric} integer $j$.
Note that the possible occupants of $[1,j]$ are $j,j+1, j+2, \dots, k(j-1)+1$.

Once we found these expressions in $n$, we can ask about the {\it limiting distribution}, that Maple can find for us.
Then we can also hope to see how it varies with $i$ and look at the {\it meta-limiting} behavior as $i$ gets larger.

Another kind of question, inspired by the beautiful work of Chan, Pak, and Panova [CPP1][CPP2], is to study the {\bf sorting probabilities}. Given two cells
$c_1$ and $c_2$,  draw a standard Young tableau {\it uniformly at random}. Who is bigger? 

The occupant of $c_1$ or the occupant of $c_2$?

Of course if the two cells are {\it related}, i.e. one of them, say $c_2$, is (weakly) below and (weakly) to the right of the other, $c_1$, i.e. in the underlying
poset $c_1<c_2$, then of course, {\it always} $T_{c_1}<T_{c_2}$, but what if they are {\bf not} related i.e. writing
$$
c_1=[i_1,j_1] \quad, \quad c_2=[i_2,j_2] \quad, \quad
$$
we have $i_1<i_2$ but $j_1>j_2$. 

The {\bf sorting probability} is defined by:
$$
SP(\lambda, c_1,c_2):=Pr(T_{c_1} > T_{c_2})-Pr(T_{c_2} > T_{c_1}) = 2\,Pr(T_{c_1} > T_{c_2})-1 \quad ,
$$
where $T$ is a random standard Young tableau of shape $\lambda$.

In particular, following Chan-Pak-Panova, we are interested in the {\bf minimal} (absolute value) of the sorting probabilities, over all pairs of cells, as the shapes get larger.

Going back to the shape $\lambda=(2,2,1)$, we see that for the first two tableaux the occupant of $(1,2)$ is less than the occupant of $(2,1)$, 
while for the last three ones it is the reverse. Hence the probability of
$T_{1,2}<T_{2,1}$ is $\frac{2}{5}$, and so $SP(221,[1,2],[2,1])=\frac{3}{5} -\frac{2}{5}=\frac{1}{5}$. The {\bf minimal sorting probability} for that shape is also $\frac{1}{5}$ (check!)

{\bf Simulation}

One way to answer these questions, {\it approximately}, is via {\it simulation}. The amazing Greene-Nijenhuis-Wilf [GrNW] algorithm (that also lead to a beautiful probabilistic proof of the {\it hook lenghth formula})
inputs a shape, $\lambda$, and outputs, {\it uniformly at random},
one of the $f^{\lambda}$ standard Young tableaux of that shape. By sampling many of them, we can get approximations to the quantities of interest.

{\eightrm
Procedure {\tt GNW(L)} in our Maple package {\tt SYT.txt} implements the Greene-Nijenhuis-Wilf algorithm. For example, try

{\tt GNW([4,3,2]);} \quad,

in order to get, uniformly at random, one of the $168$ standard Young tableaux of shape $(4,3,2)$.

To get approximations for the probability generating function, using the variable {\tt x}, of the distribution of the occupants of cell {\tt c} in a random Young tableau of shape {\tt L}, by sampling {\tt K} random tableaux, type

{\tt SiOcGF(L,c,x,K);}

For example for the shape $(4,4,4)$ and the cell $[2,2]$, with $10000$ tries, type:
{\tt SiOcGF([4,4,4],[2,2],x,10000);}

getting something like (of course it is slightly different each time)

$$
.1090000000\,x^7+.2875000000\,x^6+.3639000000\,x^5+.2396000000\,x^4 \quad.
$$

To get approximations for the sorting probability of cell {\tt c1} vs. cell {\tt c2}, in the shape {\tt L}, by sampling $K$ tableaux: enter:

{\tt SiPr(L,c1,c2,K);} \quad .

For example, 

{\tt SiPr([3,3,3],[1,2],[2,1],10000);}

would give something like {\tt 0.010400000}. Of course, in this particular case the exact answer is obviously zero, by symmetry, so getting something close to $0$ is a good {\it sanity check}.
}

{\bf Symbol Crunching in order to find The Probability Distribution of the Occupants of a Specific Cell in a Symbolic Shape}

For the sake of exposition, we will mostly be concerned with standard Young tableaux of rectangular shape. Fix the number of rows $k$, and consider the shape
$$
(n,n, \dots, n) \quad,
$$
where $n$ is repeated $k$ times. More generally all our algorithms carry over to the general {\it symbolic} shape
$$
(n_1,n_2, \dots, n_k) \quad,
$$
where $n_1 \geq n_2 \geq n_k \geq 0$, and they are all left symbolic. Using Young-Frobenius we get an {\bf explicit} formula for their total number
in the form of a certain {\bf rational function}, of $(n_1, \dots, n_k)$ times the {\bf multinomial coefficient}
$$
\frac{(n_1 + \dots + n_k)!}{n_1! \cdots n_k!} \quad .
$$
In the special case of a $k \times n$ rectangular shape:  (i.e. $(n,\dots, n)$, where $n$ is repeated $k$ times), it is a certain rational function of $n$ 
(namely $\frac{(k-1)!}{(n+1)_1 (n+1)_2 \cdots (n+1)_{k-1}}$, where, as usual $(x)_m:=x(x+1) \cdots (x+m-1)$) times
$$
\frac{(n\,k)!}{n!^k} \quad .
$$

Now also fix a specific ({\it numeric}) cell, $c=[i,j]$, and
a specific (numeric) integer $r$. We want an {\bf explicit} formula, as a {\bf rational function} of $n$, for the 
probability that, when you draw (say using the GNW algorithm) {\it uniformly at random}, one of the $\frac{(n\,k)!(k-1)!}{n!\cdots (n+k-1)!}$  standard Young tableaux,
of shape $n^k$, that the occupant of the cell $c=[i,j]$  is the integer $r$, in symbols: 

$Pr(T_{ij}=r)$.

We will soon explain how to do it, but you are welcome to try it out first using our Maple package {\tt SYT.txt}. Let's give a few examples.

$\bullet$ {\eightrm To get the  explicit expression for the probability that the occupant of cell $[1,3]$ in  a random standard Young tableau
of shape $(n,n,n)$  happens to be $7$, type:

{\tt OcCs([n,n,n],7,[1,3]);} \quad,

getting

$$
\frac{5 n \left(n +1\right)^{2} \left(n +2\right)}{9 \left(3 n -1\right) \left(3 n -2\right) \left(3 n -4\right) \left(3 n -5\right)} \quad .
$$

$\bullet$ For a more complicated example, to get the expression for the probability that the occupant of cell $[3,3]$ in  a random standard Young tableau of shape $(n,n,n)$ happens to be $13$, type:

{\tt OcCs([n,n,n],13,[3,3]);}

getting
$$
\frac{110 n^{2} \left(-1+n \right) \left(n +1\right)^{2} \left(2+n \right) \left(233 n^{2}-1933 n +3984\right)}{81 \left(-1+3 n \right) \left(-2+3 n \right) \left(-4+3 n \right) \left(-5+3 n \right) \left(-7+3 n \right) \left(-8+3 n \right) \left(-10+3 n \right) \left(-11+3 n \right)} \quad .
$$

$\bullet$ For yet another example, regarding the three-rowed shape $(n_1,n_2,n_3)$, to get the rational function (in $n_1,n_2,n_3)$ for the probability that cell $[1,2]$ would be occupied by $3$, type

{\tt OcCs([n[1],n[2],n[3]],3,[1,2]);} \quad ,

getting
$$
\frac{n_{1}^{2} n_{2}+n_{1}^{2} n_{3}+n_{1} n_{2}^{2}+2 n_{1} n_{2} n_{3}+n_{1} n_{3}^{2}+n_{2}^{2} n_{3}+n_{2} n_{3}^{2}-n_{1} n_{2}-n_{1} n_{3}+n_{2}^{2}-n_{3} n_{2}+2 n_{3}^{2}-2 n_{2}-6 n_{3}}{\left(n_{1}-2+n_{2}
+n_{3}\right) \left(n_{1}-1+n_{2}+n_{3}\right) \left(n_{1}+n_{2}+n_{3}\right)} \quad .
$$

$\bullet$
If the cell is at the first row, $c=[1,j]$, for some $j>1$,  then there are only finitely many possible occupants $r$, namely $r=j, j+1, \dots, k(j-1)+1$, and to get the probability generating function, using the variable $x$, type

{\tt OcGFs1(L,j,x);} \quad .

For example, entering:
{\tt OcGFs1([n,n,n],2,x);}

gives you
$$
\frac{2 \left(-1+n \right) x^{2}}{-1+3 n}+\frac{8 \left(-1+n \right) \left(n +1\right) x^{3}}{3 \left(-1+3 n \right) \left(-2+3 n \right)}+\frac{\left(n +1\right) \left(2+n \right) x^{4}}{3 \left(-1+3 n \right) \left(-2+3 n \right)} \quad,
$$

meaning that the cell $[1,2]$ in a standard Young tableau of shape $(n,n,n)$ is occupied by either $2$, $3$, or $4$, with respective probabilities of
$\frac{2(-1+n)}{-1+3 n}$,  $\frac{8 \left(-1+n \right) \left(n +1\right)}{3 \left(-1+3 n \right) \left(-2+3 n \right)}$, and $\frac{\left(n +1\right) \left(2+n \right)}{3 \left(-1+3 n \right) \left(-2+3 n \right)}$.

To get the limiting distribution as $n \rightarrow \infty$, as well as the expectation, variance, and the first few moments up to the $K$-th, try
{\tt OcGFs1L(L,n,i,x,K);} \quad .

}

{\bf How does Maple  perform these amazing calculations? In other words how does it compute $Pr(T_{i,j}=r)$ for a  random standard Young tableau of a symbolic shape?}

Given a (symbolic, or numeric) shape $\lambda$, a cell $c=[i,j]$, and an integer $r$, how can it happen that $T_{ij}=r$? The cells occupied by $\{1,2,\dots, r\}$ form a certain standard Young tableau with $r$ cells,
that is a certain subshape, that must contain the cell $c=[i,j]$, that must be a corner, of course. So let's ask our beloved computer to find all the  shapes with $r$ cells that contain the cell $[i,j]$ as a corner,
or equivalently the set of partitions, $\nu$, of $r$ with at least $i$ rows such that $\nu_i=j$. 
Let's call this (finite) set $S([i,j],r)$.

Let, as usual, $f^{\lambda/\nu}$ denote the number of standard Young tableau of {\bf skew-shape} $\lambda/\nu$ (recall that these are shapes where $\nu$ is a subshape of $\lambda$,
and the cells of $\nu$ are removed). Then our quantity of interest is
$$
\sum_{\nu \in S([i,j];r)} f^{\nu'} f^{\lambda/\nu} \quad ,
$$
where $\nu'$ is the shape $\nu$ with the cell $[i,j]$ removed.

The {\bf number} $f^{\nu'}$ is easily computed using the Young-Frobenius formula. How do we compute the (symbolic) {\bf expression} $f^{\lambda/\nu}$? 

Recall that standard Young tableaux of shape $\lambda=(\lambda_1, \dots, \lambda_k)$
are in bijection with {\bf walks} from the {\bf origin} to the point $\lambda$ in the $k$-dimensional discrete lattice $\N^k$, that always must stay in the region
$$
x_1 \geq x_2 \geq \dots \geq x_k  \geq 0 \quad .
$$

Similarly, standard Young tableaux of {\it skew-shape} $\lambda/\nu$ are in bijection with such `sub-diagonal' walks from $\nu$ to $\lambda$. Now following the ideas of  Andr\'e [Z] (see also [GeZ]), put
mirrors on the hyperplanes
$$
x_1 -x_2= -1 \quad, \quad
x_2 -x_3= -1 \quad, \quad
\dots \quad, \quad
x_{k-1} -x_k= -1 \quad,
$$
and look at the set of $k!$ images of the point $\nu$ under the action of the group generated by these $k-1$ reflections.
As is well-known (and fairly easy to see), the underlying group  is the symmetric group $S_k$, and the sign is $1$ or $-1$ according to whether the number of inversions is even or odd.
Calling the set of images $IMAGE(\nu)$, we have:
$$
f^{\lambda/\nu} = \sum_{\mu \in IMAGE(\nu)} \pm W(\mu,\lambda) \quad,
$$
where $W(\mu,\lambda)$ is the  {\bf number} of walks in the lattice from $\mu$ to $\lambda$ given by the {\bf multinomial coefficient}
$$
\frac{(\lambda_1+ \dots + \lambda_k - \mu_1 - \dots -\mu_k)!}{(\lambda_1-\mu_1)! (\lambda_2-\mu_2)! \cdots (\lambda_k-\mu_k)!} \quad .
$$
But since we are interested in probabilities, we can divide everything by $f^\lambda$ and stay in the realm of rational functions.

{\eightrm This is implemented in procedures {\tt Swee(L,M)}}.

{\bf Computing the Sorting Probabilities for symbolic shape and any two cells where one of them is at the first row}

{\eightrm

The {\bf numeric} procedure {\tt Pr(L,c1,c2)}, for a random standard Young tableau of shape {\tt L},
manually finds the sorting probability of cell {\tt c1} vs. cell {\tt c2}, and the numeric procedure {\tt MinPr(L)} finds
the minimal sorting probabilities among all pair of cells, followed by the `champions'. For example, if you type

{\tt MinPr([10,4,3]);} \quad,

you would get

{\tt 1/273, $\{$[[1, 5], [3, 1]]$\}$},

meaning that the minimum (absolute value) of the sorting probabilities among all the ${{17} \choose {2}}=136$ pairs of cells is $\frac{1}{127}$ and it is achieved with the pair of cells
$[1,5]$ and $[3,1]$. But we want to do things {\bf symbolically}. Alas, things get complicated if neither cells are at the first row.

But we can, exactly, and {\bf symbolically}, compute a closed-form expression, as a rational function of the symbols, of the sorting probabilities
between any cell $[1,j]$ on the first row and any cell below it (to the left, of course, or else the sorting probability is trivially $-1$).

This is implemented in procedure {\tt PrS(L,j,c2)}, where {\tt L} is the symbolic shape and {\tt c2} is the cell below the first row that we
compare it to {\tt [1,j]}. For example, to get the sorting probability of cell $[1,3]$ vs. the cell $[2,2]$, for the shape $(n,n,n)$, type

{\tt PrS([n,n,n],3,[2,2]);} getting
$$
-\frac{\left(17 n -4\right) \left(n -3\right)}{3 \left(3 n-1 \right) \left(3 n-4 \right)} \quad .
$$
}

\vfill\eject

{\bf How Does Maple find The Symbolic Sorting Probabilities?}

How do we do it? Look at all the possible occupants of cell $c_1=[1,j]$ (there are finitely many of them). Suppose it happens to be $r$. 
How can it be larger than the occupant of cell $c_2=[m_1,m_2]$?. We find  the (finite) set of shapes with $r$ cells that include $c_1=[1,j]$, and in addition
it is a corner. 
In other words all the shapes $\nu$ with $r$ cells such that $\nu_1=j$, $\nu$  has at least $m_1$ rows, and $m_2 \leq \nu_{m_1}$.

As before add-up $f^{\nu'}$ times $f^{\lambda/\nu}$, and then add-them-up for all possible legal occupants of $[1,j]$. Getting a nice (or not so nice, but sill explicit)
expression for $Pr(T_{1,j}>T_{m_1,m_2})$, and hence for the sorting probability  $2\,Pr(T_{1,j}>T_{m_1,m_2})-1$\quad.

Of course, we always divide by $f^{\lambda}$ (but this is already built-in in all our {\it macros}).

{\bf A one-line proof that the Minimal sorting probabilities for the Catalan Poset is $O(\frac{1}{n})$}

In a deep and beautiful work [CPP2], the authors proved that the minimal sorting probability of the Young lattice, as the shapes get larger, tends to $0$. In the more
specific paper [CPP1], they proved, by an ingenious and delicate asymptotic analysis, that for the two-rowed case, $[n,n]$, (what they call the {\it Catalan poset}),
it is $O(\frac{1}{n^{\frac{5}{4}}})$. But using our Maple package, we can get, {\it without human effort}, a (rigorous!) proof that it is at least $O(\frac{1}{n})$.

Indeed, entering in our Maple package {\tt SYT.txt}, the command :

{\tt PrS([n,n],3,[2,1]);} \quad ,

we get, in one {\it nano-second} :

$$
\frac{3}{2 n-1} \quad .
$$ 

So we have the following computer-generated proposition (that admittedly could have been done by humans only using paper and pencil).

{\bf Proposition}: The sorting probability of the cell $[1,3]$ and the cell $[2,1]$ in a random standard Young tableau of shape $(n,n)$ is
$$
\frac{3}{2 n-1} \, = \,\frac{3}{2} \cdot \frac{1}{n} +\frac{3}{4} \cdot \frac{1}{n^2}+\frac{3}{8} \cdot \frac{1}{n^3}+O( \frac{1}{n^4}) \quad .
$$
Hence the minimal sorting probability of the Catalan lattice is $O(\frac{1}{n})$.

Similarly, typing

{\tt PrS([n,n],5,[2,2]);}

gives the following proposition.

{\bf Proposition}: The sorting probability of the cell $[1,5]$ and the cell $[2,2]$ in a random standard Young tableau of shape $(n,n)$ is
$$
\frac{45 n^{2} - 135 n +30}{2 \left(2 n-5 \right) \left(2 n-1 \right) \left(2 n -3\right)}
\, = \,
\frac{45}{16} \cdot \frac{1}{n}+\frac{135}{32} \cdot \frac{1}{n^2} +\frac{75}{16} \cdot \frac{1}{n^3} +O(\frac{1}{n^4})
$$
Hence, again, the minimal sorting probability of the Catalan poset is $O(\frac{1}{n})$.

Procedure {\tt FindZero(L,n,K)} searches for all pairs of cells $c_1=[1,j],c_2=[m_1,m_2]$ where $c_1$ is in the first row, and $j,m_2 \leq K$, such that the sorting probability tends to $0$ (and hence is, of course $O(1/n)$).
Alas, except for the above two pairs (for the Catalan poset), none exists for $K=100$. Note that here we really lucked out, since the pairs $\{[1,3],[2,1]\}$ and $\{[1,5],[2,2]\}$ are
{\it numeric} (and small), and give upper bound for the minimal sorting probability. In order to get to the {\it true} minimum,  {\bf both} $c_1$ and $c_2$ must be {\it symbolic} 
(that what was essentially done in [CPP1] and [CPP2] with great human effort).

{\bf The special case of the Catalan poset ($2$-rowed tableaux})

For the Catalan case things can get much faster (as noticed in [CPP1]) and the procedures implementing this can be found by typing {\tt ezraD();} \quad.

{\tt Anij(n,i,j)} is a faster version of PrS([n,n],i,[2,j]). It turns out that in this case we can get closed-form expressions, for the occupancy distribution of
an arbitrary cell $[1,i]$ at the first row of a standard Young tableau of shape $(n,n)$ for {\bf symbolic} $i$, that entail, in turn,
closed-form expressions for the {\it limiting distribution} as $n$ goes to $\infty$, and {\bf surprise!} we can get explicit expressions
for the average, variance, and higher moments for that limiting distribution for symbolic $i$, and even the  {\it meta-limiting} behavior, as $i$ goes to infinity.

We have

{\bf Proposition}: The expectation of the occupant of cell $[1,i]$ in a random standard Young tableau of shape $(n,n)$, as $n$ goes to infinity is
$$
2 i +2-\frac{2 \cdot 4^{-i} \left(1+2 i \right)!}{i !^{2}} \quad,
$$
confirming Richard Stanley's observation mentioned in [CPP1], Eq. $(5.1)$ The asymptotics is
$$
2 i +2 -\frac{4}{\sqrt{\pi}} i^{1/2}
-\frac{3}{2 \sqrt{\pi}} i^{-1/2} 
+\frac{7}{32 \sqrt{\pi}} i^{-3/2} 
- \frac{9}{256 \sqrt{\pi}}\, i^{-5/2}
+ O(i^{-7/2}) \quad
$$

The variance is
$$
-\frac{4 \cdot 16^{-i} \left(1+2 i \right)!^{2}}{i !^{4}}-\frac{2 \cdot 4^{-i} \left(1+2 i \right)!}{i !^{2}}+6 i +6 \quad.
$$

The limiting (as $i$ goes to infinity) {\it skewness} is $\frac{2 \left(5 \pi -16\right) \sqrt{2}}{\left(3 \pi -8\right)^{\frac{3}{2}}}=-0.4856928234\dots $

The limiting (as $i$ goes to infinity) {\it kurtosis} is $\frac{15 \pi^{2}+16 \pi -192}{\left(3 \pi -8\right)^{2}}= 3.108163850\dots$

The limiting (as $i$ goes to infinity) {\it scaled-fifth-moment} is $\frac{2 \left(51 \pi^{2}-80 \pi -256\right) \sqrt{2}}{\left(3 \pi -8\right)^{\frac{5}{2}}}= -4.642979574\dots$

The limiting (as $i$ goes to infinity) {\it scaled-sixth-moment} is  $\frac{105 \pi^{3}+648 \pi^{2}-2240 \pi -2560}{\left(3 \pi -8\right)^{3}}=18.66866547\dots$

For more details see the output file

{\tt https://sites.math.rutgers.edu/\~{}zeilberg/tokhniot/oSYT9new.txt} \quad .

We note that, by miracle, the (limiting, as $n$ goes to infinity) average, variance, and {\it any} higher moment, happened to be {\it gosperable} so Maple
is able to evaluate them in closed-form using the Maple command {\tt sum}. We doubt whether this will happen for more rows, but we did not try.

We believe that the Maple package {\tt SYT.txt} can be used to explore further and possibly suggest improvements to the already very impressive work in [CPP1] and [CPP2].

{\bf Sample Output}

The web-page of this article

{\tt https://sites.math.rutgers.edu/\~{}zeilberg/mamarim/mamarimhtml/syt.html} \quad ,

contains lots of output file. Let's just mention some highlights.

$\bullet$ If you want to see a computer-generated article with lots of explicit expressions (as rational functions of $n$) 
for the probability distribution of the occupant of cell $[1,i]$ in a (uniformly-at) random-generated Young tableau of rectangular shape with $2$ rows and $n$ columns (i.e. of shape $[n,n]$) 
for all $i$ between $2$ and $40$, see the output file

{\tt https://sites.math.rutgers.edu/\~{}zeilberg/tokhniot/oSYT1.txt} \quad .

If you want to see an abbreviated version, with only the {\it limiting distribution} as $n$ goes to infinity, but all the way to $i=60$, as well as the
expectation, variance, skewness, and kurtosis, see

{\tt https://sites.math.rutgers.edu/\~{}zeilberg/tokhniot/oSYT1L.txt} \quad .

If you want to see the plot of the occupancy distribution of the cell $[1,40]$ as $n$ goes to $\infty$, look here:

{\tt https://sites.math.rutgers.edu/\~{}zeilberg/tokhniot/oSYTpic2.html} \quad .

$\bullet$ For the analogous output for $3$-rowed rectangular shapes, (i.e. $(n,n,n)$), see, respectively

{\tt https://sites.math.rutgers.edu/\~{}zeilberg/tokhniot/oSYT2.txt} \quad,

{\tt https://sites.math.rutgers.edu/\~{}zeilberg/tokhniot/oSYT2L.txt} \quad,

{\tt https://sites.math.rutgers.edu/\~{}zeilberg/tokhniot/oSYTpic3.html} \quad .

For the analogous information for up to $8$ rows (but with less data) see the above-mentioned front of this article.

$\bullet$ For testing the amazing Greene-Nijenhuis-Wilf algorithm vs. the exact results, see the output file

{\tt https://sites.math.rutgers.edu/\~{}zeilberg/tokhniot/oSYTsi1.txt} \quad .

To get lots of explicit expressions for sorting probabilities refer to the above web-page.

{\bf References}

[CPP1] Swee Hong Chan, Igor Pak, and Greta Panova, {\it Sorting probabilities of Catalan Posets},  Advances in Applied Mathematics {\bf 129} (2021), article 102221, 13 pp.\hfill\break
{\tt https://www.math.ucla.edu/\~{}pak/papers/Cat-sort15.pdf} .

[CPP2] Swee Hong Chan, Igor Pak, and Greta Panova, {\it  Sorting probability for large Young diagrams}, Discrete Analysis, Paper 2021:24, \hfill\break
{\tt https://www.math.ucla.edu/\~{}pak/papers/Sorting20.pdf} \quad .

[F] William Fulton, {\it ``Young Tableaux, with Applications to Representation Theory and Geometry''}, Cambridge University Press, 1997.

[K] Donald E. Knuth,  {\it ``The Art of Computer Programming, Vol. III: Sorting and Searching (2nd ed.)''}, Addison-Wesley, (1973), (section 5.1.4).

[GeZ] Ira Gessel and Doron Zeilberger, {\it Random walk in a Weyl chamber}, Proc. Amer. Math. Soc. {\bf 115} (1992), 27-31. \hfill\break
{\tt https://sites.math.rutgers.edu/\~{}zeilberg/mamarim/mamarimPDF/weyl.pdf} \quad .

[GrNW] Curtis Greene, Albert Nijenhuis, and Herbert Wilf, {\it A probabilistic proof of a formula for the number of Young Tableaux of a given shape}, Adv. in Math. {\bf 31} (1979), 104-109. \hfill\break
{\tt https://www2.math.upenn.edu/\~{}wilf/website/Probabilistic\%20proof.pdf} \quad .

[Wi] Wikipedia, {\it Young tableaux}, {\tt https://en.wikipedia.org/wiki/Young$\_$tableau} \quad .

[Z] Doron Zeilberger, {\it Andr\'e's reflection  proof generalized to the
many-candidate ballot problem},  Discrete Mathematics {\bf 44} (1983), 325-326. \hfill\break
{\tt https://sites.math.rutgers.edu/\~{}zeilberg/mamarimY/Andre1983.pdf} \quad.  \hfill\break
{\eightrm [Christian Krattenthaler noticed a very long time ago that on the second page (p. 326), last word on line 3: first $\rightarrow$ last]} \quad .

\vfill\eject

\bigskip
\hrule
\bigskip
Shalosh B. Ekhad and Doron Zeilberger, Department of Mathematics, Rutgers University (New Brunswick), Hill Center-Busch Campus, 110 Frelinghuysen
Rd., Piscataway, NJ 08854-8019, USA. \hfill\break
Email: {\tt  ShaloshBEkhad at gmail dot com} \quad, \quad {\tt DoronZeil at gmail dot com}   \quad .

Written: {\bf March 29, 2023}.

\end